\documentclass[a4paper,12pt]{article}
\usepackage{hyperref}
\usepackage{amssymb}
\usepackage{amsmath}
\usepackage[latin1]{inputenc} 
\newtheorem{defi}{Definition}[section]
\newtheorem{prop}[defi]{Proposition}
\newtheorem{thm}[defi]{Theorem}
\newtheorem{lem}[defi]{Lemma}
\def\R{{\mathbb R}}
\def\C{{\mathbb C}}
\def\Z{{\mathbb Z}}
\def\N{{\mathbb N}}
\def\Q{{\mathbb Q}}
\def\F{{\mathbb F}}
\def\OO{{\mathcal O}}

\def\eps{{\varepsilon}}

\def\supp{{\tx{supp}}}
\def\BB{{\mathcal B}}
\def\Ker{{\tx{Ker}}}
\def\strongT{{property $(T_{Schur},G,K)$ }}
\def\tgk{$(T_{Schur},G,K)$ }

\newcommand {\tx}[1] {\textrm{#1}}
\newcommand {\resr}[1] {\OO/\pi^{#1}\OO}

\newcommand {\esp}[1]{\underset{ #1}{\mathbb E}}
\newcommand{\qed}{\hfill $\Box$}
\newcommand{\pair}[1]{\langle{#1}\rangle}

\newcommand{\sch}[1]{\tx{Schur}_{#1}}

\begin{document}

\title{About the difficulty to prove the Baum Connes conjecture without coefficient for a non-cocompact lattice in $Sp_4$ in a local field}
\author{
Benben Liao\footnote{
Institut de Mathématiques de Jussieu - Paris Rive Gauche
(Université Paris Diderot - Paris 7)
}} 
\date{\today}
\maketitle
\abstract{We introduce property \tgk 
and prove it
for some non-cocompact lattice in $Sp_4$ in a local field of finite characteristic. We show that 
property \tgk for a non-cocompact lattice $\Gamma$ in a higher rank almost simple algebraic group in a local field
is an obstacle to proving 
the Baum Connes conjecture without coefficient for $\Gamma$
with known methods,
 and 
this is stronger than the well-known fact that
$\Gamma$ does not have the property of rapid decay (property (RD)).
It is the first example (as announced in \cite{laff-orsay})
for which all known methods for proving the Baum Connes conjecture without coefficient fail.
}

\section{Introduction}

N.~Higson, V.~Lafforgue, and
G.~Skandalis constructed 
counterexamples to the Baum Connes conjecture for discrete group actions on commutative $C^*$ algebras 
 using Gromov's groups which do not uniformly embed into Hilbert space \cite{counterex}. 
V.~Lafforgue introduced strong Banach property $(T)$ \cite{laf-duke,laf-jta}, proved it for $SL_3(\Q_p)$, and constructed the
first example of expander graphs which do not embed uniformly into any Banach space of non-trivial type.
Other examples of expander graphs 
non-embeddable in Banach spaces of non-trivial type or Banach spaces of weaker properties 
have been found 
\cite{mendel-naor,liao,delasalle}.
In \cite{laff-ncg}, V.~Lafforgue introduced property $(T_{Schur})$ 
(which is stronger than strong property $(T)$ \cite{laf-duke}) 
and proved that it is an obstacle 
to proving Baum Connes conjecture for $SL_3(\Q_p)$ with commutative coefficient containing $C_0(SL_3(\Q_p))$ with known methods.
In this article, we introduce 
property \tgk in Definition~\ref{strongT}
as an analogue of property $(T_{Schur})$, prove it for the non-cocompact lattice $\Gamma=Sp_4(\F_q[\pi^{-1}])$ of $Sp_4(\F_q((\pi)))$ in Theorem~\ref{main} (which is the main result of this article), and show
that it is an obstacle to proving Baum Connes conjecture 
\textit{without coefficients}
for the lattice $\Gamma$ with known methods,
which is stronger than the well-known fact that
$\Gamma$ does not have the property of rapid decay (property (RD)).
It is the first example for which all known methods for proving Baum Connes conjecture 
\textit{without coefficients}
fail.

We begin with some notations, and then state the main result of this article.

Let $G$ be non-compact locally compact topological group, $K\subsetneq G$ a compact subgroup.
Let $H\subseteq G$ be a non-compact closed subgroup. 
Let $\ell:G\to\R_{\geq 0}$ be a continuous length function on $G.$
Denote $B_n$ the ball of radius $n$ in $G.$

For any continuous function $c\in C(G)$ on $G,$
we introduce the following notation for the norm of the Schur product by $c$ on the subspace 
$C(H\cap B_n)=\{f\in C_c(H),\supp(f)\subset B_n\}$
of functions on $H$ with supports in $B_n,$

$$\|\sch{c}|_{C(H\cap B_n)}\|=\|\sch{c}|_{C(H\cap B_n)}\|_{\mathcal L(C^*_r(H))}$$
$$=\sup\{\|\sch{c}f\|_{C^*_r(H)},f\in C_c(H),\supp(f)\subset B_{n},\|f\|_{C^*_r(H)}\leq 1\},$$
where $\sch{c}f\in C_c(H)$ denotes the Schur product 
$$\sch{c}f(h)=c(h)f(h),\forall h\in H.$$

\begin{defi} \label{strongT}
We say that $H$ has \strongT 
if for any continuous length function $\ell:G\to\R_{\geq 0},$
there exists $s_0>0$ such that $\forall s\in[0,s_0)$
there exists a continous function $\phi\in C_0(G)$ vanishing at infinity,
such that 
$\forall C>0$ and 
for any family of $K$-biinvariant functions $c\in C(G)$ with the following uniform Schur condition
\begin{equation}\label{condition-schur}
\|\sch{c}|_{C(H\cap B_n)}\|_{\mathcal L(C_r^*(H))}\leq Ce^{sn},\forall n\in\N,
\end{equation} 
there exists a limit $c_\infty \in\C$ to which $c$ tends uniformly rapidly 
$$|c(g)-c_\infty|\leq C\phi(g),\forall g\in G.$$  
\end{defi}

Let $\F_q$ be a finite field of characteristic different from $2$ with
cardinality $q.$ 
Let $G$ be 
$Sp_4(\F_q((\pi)))$ over the local field $\F_q((\pi)),$
$K=Sp_4(\F_q[[\pi]])$ a maximal compact subgroup of $G.$
Let $\Gamma$ be the non-cocompact lattice $Sp_4(\F_q[\pi^{-1}])$ in $G.$
Let $H\subsetneq\Gamma$ be the unipotent subgroup consisting of elements of the form 
$$\begin{pmatrix}1 &*&*&*\\
 & 1 &*&*\\
 & & 1 &*\\
 & & & 1
\end{pmatrix}\in \Gamma.$$

The following is the main result of this article.

\begin{thm}\label{main}
The unipotent group $H$ has \strongT as in Definition~\ref{strongT}.
As a consequence, 
the 
non-cocompact
lattice $\Gamma=Sp_4(\F_q[\pi^{-1}])\subsetneq G=Sp_4(\F_q((\pi)))$ also has \strongT .
\end{thm}
\noindent
{\bf Remark.} That $\Gamma$ has \strongT follows from that $H$ has property 
\mbox{\tgk.} 
Indeed, 
it is clear from Definition~\ref{strongT} that for any two discrete subgroups $H\subsetneq H'\subseteq G,$
if $H$ has \strongT, then $H'$ also has 
\mbox{\strongT.}

$H$ is an amenable group and thus satisfies the Baum Connes conjecture. 
Theorem~\ref{main} says that
the constant function 
$1$ on $H$ cannot be deformed among $K$-biinvariant functions on $G$ satisfying the Schur type condition \eqref{condition-schur}, whereas it does not prevent
the possibilities of deformations among other functions.

For the non-cocompact lattice $\Gamma=Sp_4(\F_q[\pi^{-1}]),$ 
\strongT
is an obstacle of known methods for
proving Baum Connes conjecture (without coefficient).
It plays the role of property $(T_{Schur})$ 
for a locally compact group $G$ (e.g. $SL_3(\Q_p)$)
relative to some open compact subgroup (e.g. $SL_3(\Z_p)$)
as in \cite{laff-ncg}
being an obstacle of known methods for proving Baum Connes conjecture for $G$ with commutative coefficient containing $C_0(G)$. 
Indeed,
if a lattice in a higher rank almost simple algebraic group has \mbox{\strongT,}
then it does not have property (RD) 
(Proposition~\ref{prop-RD}).
A.~Valette conjectured that any cocompact lattice in a simple algebraic group over a local field has property (RD), thus cocompact lattice wouldn't seem to have \strongT.
What is shown in this article is stronger:
\strongT for a higher rank lattice $\Gamma$ prevents the existence of any reasonable dense subalgebra of $C^*_r(\Gamma)$ for known methods for proving the Baum Connes conjecture (Proposition~\ref{prop-contradiction}), in particular the Jolissaint algebra for property (RD).
We recall known methods for proving Baum Connes conjecture for a lattice in a reductive group over a local field in section 2, and prove in Proposition~\ref{prop-contradiction} that the conditions of these methods are not satisfied for any lattice with \strongT
(by adapting the arguments in \cite{laff-ncg} to our situation).

We will give two proofs of Theorem~\ref{main} in section 3 and section 4 respectively. 

The first proof is more in line with the arguments in \cite{laff-preprint}, and is therefore more transparent and more checkable. What is different from \cite{laff-preprint} is the use of two families of parameters (i.e. $n=2$ in Lemma~\ref{lem-single-pair}) which yields an improved estimate $q^{-j}$ in
the second inequality of
 Proposition~\ref{estimates-worse}.

The second proof makes use of two families of functions on $H$ with exponentially small $C^*_r$ norms, and yields a slightly better constant $s_0$ in Definition~\ref{strongT}. 
The first family of functions corresponding to the abelian subgroups are already constructed and the corresponding estimate is obtained in
\cite{laff-orsay}, which we include in this article using the same arguments. 
It is constructed in this article the other family of functions corresponding to the discrete Heisenberg subgroup. 
The improved estimate in $q^{-j}$ is obtained using harmonic analysis on the Heisenberg group over the ring of polynomials $\F_q[\pi^{-1}].$

{\bf Acknowledgement}
I would like to thank Vincent Lafforgue for his useful suggestions and precise explanations of various aspects on this problem. I thank Georges Skandalis for the discussion on algebra $A_\theta.$
I thank Mikael de la Salle for several corrections to this article, especially for
pointing out a mistake in the parameterization of the group of Heisenberg type.

\section
{An obstacle to proving the Baum Connes conjecture without coefficients}

In this section,
$G$ is an arbitary reductive group over a local field, $K\subsetneq G$ a maximal compact subgroup, 
$\ell:G\to \R_{\geq 0}$ a continuous length function, and
$\Gamma \subsetneq G$ a lattice.

Adapting conditions $(\tilde D)=(D1)+(D2)+(D3)+(D4)$ in \cite{laff-ncg} to our situation, we list the following conditions $(\tilde D')=(D1)+(D2)+(D3)+(D4')$ to include all known methods to prove Baum Connes conjecture for $\Gamma$. More precisely, conditions $(D1),(D2),(D3)$ are just specialization to $A=\C$ of conditions $(D1),(D2),(D3)$ \cite{laff-ncg}, and condition $(D4')$ is a variant of $(D4)$ \cite{laff-ncg} in which we require the representations of $\Gamma$ for the homotopy come from representations of $G.$

\begin{itemize}
\item $(D1)$. For any $s>0,$ there exists $C_s>0,$ and a Banach subalgebra $\BB_s\subset C_r^*(\Gamma)$ containing $\C(\Gamma)$ as a dense subalgebra, such that $\forall n\in\N,\forall f\in\C(\Gamma)$ with $\supp(f)\subset B_n,$
$$\|f\|_{\BB_s}\leq C_se^{sn}\|f\|_{C^*_r(\Gamma)},$$
where $B_n\subset \Gamma$ denotes the ball of length $n.$
\item $(D2)$. There exist a homotopy $(E,\pi,T)\in E_{\Gamma,?}^{ban}(\C,\C[0,1])$ ($?$ indicates that there is no restriction on the norms of the representations of $\Gamma$) from $1$ to the $\gamma$ element, and $(\tilde E,\tilde \pi,\tilde T)\in E^{ban}(\BB_s,C^*_r(\Gamma)[0,1])$ with $(\tilde E^<,\tilde E^>)$ containing $(\C(\Gamma,E^<),\C(\Gamma,E^>))$ as a dense subspace, such that the imbeddings 
$$i_s:\C(\Gamma)\hookrightarrow\BB_s,i_r:\C(\Gamma,\C[0,1])\hookrightarrow C_r^*(\Gamma,\C[0,1]),$$
$$i^<:\C(\Gamma,E^<)\hookrightarrow\tilde E^<,i^>:\C(\Gamma,E^>)\hookrightarrow\tilde E^>$$ 
and 
$$\pair{\cdot,\cdot}:\tilde E^<\times \tilde E^>\to C^*_r(\Gamma)[0,1]$$ 
satisfy: $\forall f\in\C(\Gamma),\varphi\in\C(\Gamma,\C[0,1]),F_1\in\C(\Gamma,E^<),F_2\in\C(\Gamma,E^>),$
$$i^<(\varphi F_1)=i_r(\varphi)i^<(F_1),i^<(F_1f)=i^<(F_1)i_s(f),$$
$$i^>(fF_2)=i_s(f)i^>(F_2),i^>(F_2\varphi)=i^>(F_2)i_r(\varphi),$$
and
$$i_r(\pair{F_1,F_2})=\pair{i^<(F_1),i^>(F_2)}.$$
\item $(D3)$. The $\C[0,1]$ pair $(E^<,E^>)$ is in isometric duality as defined in \cite{laff-ncg} (definition 2.13), i.e. the maps 
$$E^<\to\mathcal L_{\C[0,1]}(E^>,\C[0,1]),E^>\to\mathcal L_{\C[0,1]}(E^<,\C[0,1])$$
 are isometries.
\item $(D4')$. The representation $E$ of $\Gamma$ is restriction of some representation of $G$ on which $K$ acts by isometries, and $\forall x\in E^>,\xi\in E^<,$
$$\|e_0\otimes x\|_{\tilde E^<}\leq \|x\|_{E^<,}\|e_0\otimes \xi\|_{\tilde E^>}\leq \|\xi\|_{E^<}.$$
\end{itemize} 

\begin{prop} \label{prop-contradiction}
Suppose $\Gamma$ has \strongT as in Definition~\ref{strongT}. 
 Then for any continuous length function $\ell$ on $G,$
 $G,K,\ell,\Gamma$ do not satisfy $(\tilde D')=(D1)+(D2)+(D3)+(D4').$ 
\end{prop}

The proof of Proposition~\ref{prop-contradiction} is an adaptation of the proof of Proposition~4.2 in \cite{laff-ncg} to our situation.

\begin{lem}\label{fund-cal}
Suppose $\Gamma$ satisfies $(D1)+(D2)+(D4').$ Then $\forall s>0,$ there exists $C_s>1,$ such that $\forall x\in E^>,\xi\in E^<$ both of norms $\leq 1,$ putting $c_t(\gamma)=\pair{\xi,\pi_t(\gamma)x},\forall\gamma\in\Gamma,$ we have
$\forall n\in\N,\forall t\in [0,1],$
$$\|\sch{c_t}|_{\C(\Gamma\cap B_n)}\|\leq C_s e^{sn}.$$
\end{lem}

{\bf Proof.} 
For any $f\in\C(\Gamma)$ we have the following fundamental calculation
$$\pair{e_0\otimes\xi,f(e_0\otimes x)}=\sum_{\gamma\in\Gamma}f(\gamma)\pair{\xi,\pi_t(\gamma)x}e_\gamma$$
$$=\sch{c_t}f\in\C(\Gamma).$$ By condition $(D2)$ and $(D4'),$
$$\|\sch{c_t}f\|_{C_r^*(\Gamma)}\leq \|e_0\otimes\xi\|_{\tilde E^<}\|f\|_{\BB_s}\|e_0\otimes x\|_{\tilde E^>}$$
$$\leq\|\xi\|_{E^<}\|f\|_{\BB_s}\|x\|_{E^>}\leq\|f\|_{\BB_s}.$$
When $\supp f\subset B_n,$ by condition $(D1),$
$$\|\sch{c_t}f\|_{C_r^*(\Gamma)}\leq C_se^{sn}\|f\|_{C_r^*(\Gamma)}.$$
\qed

{\bf Proof of Proposition~\ref{prop-contradiction}}
We prove it by contradiction. 
Suppose $\Gamma,G,K,\ell$ satisfy $(\tilde D')$ and $\Gamma$ has \strongT.
By Lemma~\ref{fund-cal} we see that when $\xi,x$ are both $K$-invariant,
$c_t(g)$ is a Cauchy sequence
$$|c_t(g)-c_t(g')|\leq C_s(\phi(g)+\phi(g')),\forall g,g'\in G.$$
By condition $(D3)$ we have
$$\sup_{\|\xi\|_{E^<}\leq 1,K\tx{-invariant,}t\in [0,1]}|\pair{\xi,(\pi_t(g)-\pi_t(g'))x}|$$
$$=\|\pair{\cdot,(\pi(e_Ke_g)-\pi(e_Ke_{g'}))x}\|_{\mathcal L(E^<,\C[0,1])}$$
$$=\|(\pi(e_Ke_g)-\pi(e_Ke_{g'}))x\|_{E^>},$$
as a consequence $\pi(e_Ke_ge_K)$ is also a Cauchy sequence
$$\|\pi(e_Ke_ge_K)-\pi(e_Ke_{g'}e_K)\|_{\mathcal L_{\C[0,1]}(E^<)}\leq C_s(\phi(g)+\phi(g')).$$
For the same reason 
$$\|\pi(e_Ke_ge_K)-\pi(e_Ke_{g'}e_K)\|_{\mathcal L_{\C[0,1]}(E^>)}\leq C_s(\phi(g)+\phi(g')).$$
Denote by $P$ the limit of $\pi(e_Ke_ge_K).$ 
We see that $g'kg$ tends to infinity when $g'$ tends to infinity
since $\ell(gkg')\geq \ell(g')-\ell(g^{-1}).$
Therefore, we have 
$$e_Ke_gP=\lim_{g'}\pi(e_K\int_Ke_{gkg'}dke_K)=P,$$
and 
$$P^2=\lim_g\lim_{g'}\pi(e_K\int_Ke_{gkg'}dke_K)=P.$$
Moreover, when $E_t$ is a Hilbert space and
$(\pi_t,E_t)$ is a unitary representation of $G,$ $P_t\in\mathcal L(E_t)$ is the projection onto $G$-invariant vectors
$P_t E_t=E_t^G.$
Indeed,
$\forall x\in P_tE_t,$ $\forall g\in G,$
 $$\|\pi(e_g)x-\pi(e_Ke_g)x\|^2=\|\pi(e_g)x\|^2-\|\pi(e_Ke_g)x\|^2=\|x\|^2-\|x\|^2=0,$$ 
 we have 
 $$x=\pi(e_Ke_g)x=\pi(e_g)x.$$

$P\in\mathcal L_{\C[0,1]}(E)$ and consequently $PTP\in\mathcal L_{\C[0,1]}(E).$ We denote by $\tx{Im}P$ the $\C[0,1]$-pair whose underlying Banach spaces are the images of $E^<,E^>$ under the maps $P^<,P^>.$ We have that $(\tx{Im}P,PTP)\in E^{ban}(\C,\C[0,1]).$ Indeed, 
$[e_Ke_ge_K,T]\in\mathcal K_{\C[0,1]}(E)$, as a consequence $[P,T]\in\mathcal K_{\C[0,1]}(E).$ 
Moreover,
$$P-(PTP)^2P=P(1-T^2)+P[P,T]T+PTP[P,T]\in\mathcal K_{\C[0,1]}(E),$$
which means $Id_{\tx{Im}P}-(PTP)^2\in\mathcal K_{\C[0,1]}(\tx{Im}P).$

Now $(\pi_0,E_0)$ is the trivial representation of $G$ $(E_0=\C),$ and $(\pi_1,E_1)$ is a unitary representation of $G$ without $G$-invariant vectors. $P_0T_0P_0:\C\to 0$ has index $1$ whereas $P_1T_1P_1:0\to 0$ has index $0,$ this is a contradiction and the proposition is proved.
\qed

Now let $G$ be a semisimple group over a local field, $K\subsetneq G$ a maximal compact subgroup, $\Gamma\subsetneq G$ a lattice.

Let $\ell:G\to\R_{\geq 0}$ be the $K$ biinvariant length function induced from the $G$-invariant Riemannian metric from the symmetric space or the Bruhat-Tits building associated to $G.$
Recall that when the split rank of $G$ is $\geq 2,$ $\Gamma$ has Kazhdan's property $(T)$ and thus is finitely generated. The word metric and Riemannian metric on $\Gamma$ are bi Lipchitz 
\cite{lub}.

When $\Gamma$ has property $(RD),$ it is shown in \cite{laff-ncg} that $G,K,\ell,\Gamma$ satisfy conditions $(\tilde D')$ above, and thus do not fulfil the condition in definition \ref{strongT}. We give a direct proof of this fact.

\begin{prop}\label{prop-RD}
If $\Gamma$ has \strongT as in definition \ref{strongT},
then $\Gamma$ does not have property (RD) for 
any continuous length function restricted from $G.$
In particular when the split rank of $G$ is $\geq 2,$ 
$\Gamma$ does not have property (RD) for the word length.
\end{prop}
{\bf Proof.} 
First, for any function $f\in\C(\Gamma)$ we have
\begin{equation}\label{ineq-ell2}
\|f\|_{\ell^2(\Gamma)}\leq\|f\|_{C_r^*(\Gamma)}.
\end{equation}
Indeed, for such $f,$ $$\|f\|_{\ell^2(\Gamma)}=\|fe_0\|_{\ell^2(\Gamma)}$$
$$\leq\|f\|_{C_r^*(\Gamma)}\|e_0\|_{\ell^2(\Gamma)}=\|f\|_{C_r^*(\Gamma)}.$$

Now suppose $\Gamma$ has property (RD), i.e. there exist $R>0,D\geq 0,$ such that $\forall f\in\C(\Gamma)$ with $\supp f\subset B_n,$
\begin{equation}\label{ineq-RD}
\|f\|_{C_r^*(\Gamma)}\leq Rn^D\|f\|_{\ell^2(\Gamma)}.
\end{equation}

Denote $\chi_{B_m}$ the characteristic function of $B_m$ (ball of radius $m$) for $m\in\N.$ 
For $f\in\C(\Gamma)$ 
with $\supp f\subset B_n,$ apply
\eqref{ineq-RD} to $\sch{\chi_{B_m}}f\in\C(\Gamma)$ and \eqref{ineq-ell2} to $f,$
$$\|\sch{\chi_{B_m}}f\|_{C_r^*(\Gamma)}$$
$$\leq R\min(m,n)^D \|\sch{\chi_{B_m}}f\|_{\ell^2(\Gamma)}\leq Rn^D\|f\|_{\ell^2(\Gamma)}$$
$$\leq Rn^D\|f\|_{C^*_r(\Gamma)},\forall m\in\N.$$
Namely, for any $s>0,$ there exists $C_s>0$ such that
$$\|\sch{\chi_{B_m}}|_{\C(\Gamma\cap B_n)}\|\leq C_se^{sn},\forall n\in\N.$$
Now let $s\in(0,s_0),$ by the hypothesis in the proposition, 
$$|\chi_{B_m}(\gamma)|\leq C_s\phi(\gamma),\forall m\in\N,\forall \gamma\in\Gamma,$$
which is a constadiction to the assumption that $\phi\in\C_0(\Gamma)$ is a function vanishing at infinity.\qed

\section{First proof of Theorem~\ref{main}} 
First let us be more precise on the notations.

Let by $\F_q$ a finite field of characteristic different from $2$ with
cardinality $q.$ Denote $F=\F_q((\pi))$ the local field of Laurent series in $\pi$ with coefficients in $\F_q,$ $\OO=\F_q[[\pi]]$ the ring of formal series in $\pi,$ i.e. the ring of integers of $F.$

Let $G=Sp_4(F),$ i.e. the symplectic group of $4$ by $4$ matrices $A\in M_{4}(F)$ satisfying $AJ^tA=J,$ where $^tA$ denotes the transpose of $A$ and 
$$J=\begin{pmatrix}0&0&0&1\\0&0&1&0\\0&-1&0&0\\-1&0&0&0\end{pmatrix}.$$ 
Denote $K=Sp_4(\OO),$ a maximal compact subgroup of $G.$ 

Now denote $\Gamma=Sp_4(\F_q[\pi^{-1}]).$ By the well-known reduction theory of Harish-Chandra-Borel-Behr-Harder, $\Gamma$ is lattice in $G,$ and
by Goldman's compactness criterion (see IV 1.4 in \cite{oes} in the case of characteristic $p$), it is non-cocompact.

Denote $H$ the unipotent subgroup in $\Gamma$ consisting of elements of the form 
$$\begin{pmatrix}1 &*&*&*\\
 & 1 &*&*\\
 & & 1 &*\\
 & & & 1
\end{pmatrix}\in \Gamma.$$

We define an explicit length function on $G.$
Denote by $D(i,j)$ the diagonal element $$D(i,j)=\begin{pmatrix}\pi^{-i} & & & \\
 & \pi^{-j} & & \\
 & & \pi^j & \\
 & & & \pi^i
\end{pmatrix}\in G$$ for any $(i,j)$ in the Weyl chamber $\Lambda=\{(i,j)\in\N,i\geq j\}.$
By Cartan decomposition, $\Lambda$ is in bijection with the double coset $K\backslash G/K$ via the map $(i,j)\mapsto KD(i,j)K,(i,j)\in\Lambda.$ Let $\ell:G\to\N$ be the length function on $G$ defined by $\ell(kD(i,j)k')=i+j,(i,j)\in\Lambda,k,k'\in K.$ It is clear that $\ell$ is bi-Lipschitz to the length function induced from the distance on the Bruhat-Tits building associated to $G.$ 
For any continuous length function $\ell':G\to\R_{\geq 0},$ there exists $\kappa>0$ such that $\ell'\leq\kappa(\ell+1).$

It is clear that $H$ surjects onto the double cosets $K\backslash G/K,$ since
$$\begin{pmatrix}
1&0&0&\pi^{-j}\\
&1&\pi^{-i}&0\\
 &&1&0\\
 & & &1
\end{pmatrix}\in KD(i,j)K.$$

The following theorem is a more precise statement of Theorem~\ref{main} in the introduction.

\begin{thm}\label{obstruction}
Let $s_0=\frac{\log q}{6}.$ There exists a constant $C_q>0$ depending only on $q,$ such that the following holds. For any $s\in[0,s_0)$ and $C>0,$ for any $K$-biinvariant function $c\in C(G)$ with the following Schur condition
\begin{equation}\label{schur}
\|\sch{c}|_{\C(H\cap B_n)}\|
\leq Ce^{sn}, \forall n\in\N,
\end{equation} 
there exists a limit $c_\infty\in\C$ to which $c$ tends exponentially fast
$$|c(g)-c_\infty|\leq CC_qe^{-\ell(g)(s_0-s)/4},\forall g\in G.$$
\end{thm}

{\bf Proof of Theorem.~\ref{main} by
Theorem.~\ref{obstruction}.}
Let $\ell'$ be any length function on $G.$ There exists $\kappa>0$ such that 
$\forall g\in G, \ell'(g)\leq \kappa(\ell(g)+1).$ With $s_0'=s_0/\kappa,s'\in[0,s_0')$ and
$\phi'(g)=C_qe^{s\kappa}e^{-\ell(g)(s_0-s\kappa)/4},\forall g\in G,$
Theorem.~\ref{main} is proved.
\qed

\begin{prop}\label{estimates-worse}
For any $K$ biinvariant function $c\in C(G),$ we have
\begin{equation}
|c(D(i,j))-c(D(i,j+1))|\leq 2q^{-(i-j)/2}\|\sch{c}|_{\C(H\cap B_{2i})}\|,
\end{equation}
for any $(i,j)\in\Lambda$ with $i\geq 1$ and
\begin{equation}
|c(D(i,j))-c(D(i+1,j-1))|\leq 2q^{3}\cdot q^{-j}\|\sch{c}|_{\C(H\cap B_{i+j})}\|,
\end{equation}
for any $(i,j)\in\Lambda$ with $j\geq 3.$
\end{prop}

{\bf Proof of Theorem~\ref{obstruction} by Proposition~\ref{estimates-worse}:}

By hypothesis
$$\|\sch{c}|_{\C(H\cap B_n)}\|_{C_r^*(H)}\leq Ce^{sn},\forall n\in\N,$$
and thus the two inequalities imply respectively
\begin{equation}
|c(D(i,j))-c(D(i,j+1))|\leq 2q^{-(i-j)/2}Ce^{2si},
\end{equation}
\begin{equation}
|c(D(i,j))-c(D(i+1,j-1))|\leq 2q^{3}\cdot q^{-j}Ce^{s(i+j)}.
\end{equation}
And in particular
$$|c(D(3j,j))-c(D(3j+3,j+1))|\leq CC_qe^{-(\log q-6s)j}.$$
By moving along the line $i=3j$ in the Weyl chamber as illustrated below, we have
$$|c(D(3j,j))-c_\infty|\leq CC_qe^{-(\log q-6s)j}.$$
When $i\geq 3j,$
$$|c(D(i,j))-c_\infty)|\leq CC_qe^{-(\log q -6s)i/3}\leq CC_qe^{-(\log q -6s)(i+j)/4}.$$
When $i\leq 3j,$
$$|c(D(i,j))-c_\infty)|\leq CC_qe^{-(\log q -6s)(i+j)/4}.$$

\vskip 10mm

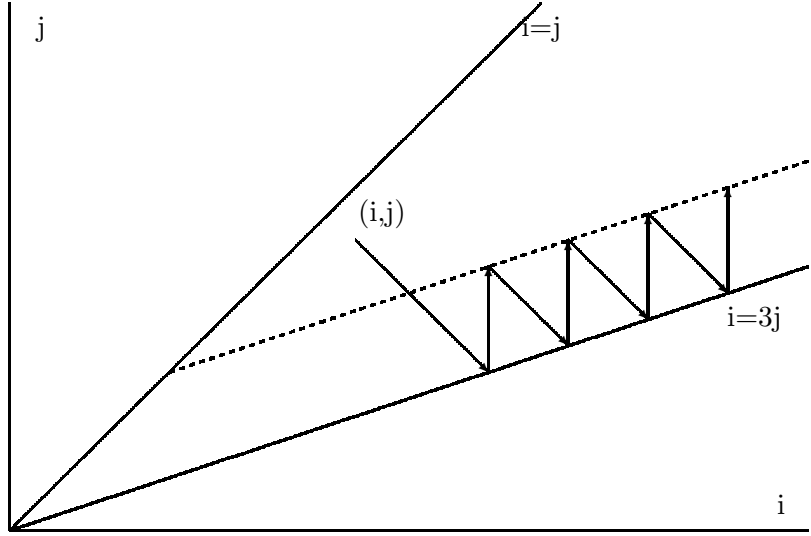
 \def\JPicScale{0.7}
\ifx\JPicScale\undefined\def\JPicScale{1}\fi
\unitlength \JPicScale mm
\begin{picture}(160,120)(0,0)
\linethickness{0.3mm}
\put(0,10){\line(1,0){150}}
\linethickness{0.3mm}
\multiput(0,10)(0.12,0.12){833}{\line(1,0){0.12}}
\put(80,85){\makebox(0,0)[cc]{}}

\linethickness{0.3mm}
\put(0,10){\line(0,1){100}}

\put(70,70){\makebox(0,0)[cc]{(i,j)}}

\put(145,15){\makebox(0,0)[cc]{i}}

\linethickness{0.3mm}
\multiput(0,10)(0.36,0.12){417}{\line(1,0){0.36}}
\linethickness{0.3mm}
\multiput(30,40)(1.89,0.63){64}{\multiput(0,0)(0.31,0.1){3}{\line(1,0){0.31}}}
\put(35,-5){\makebox(0,0)[cc]{}}


\put(70,10){\makebox(0,0)[cc]{}}

\put(6,105){\makebox(0,0)[cc]{j}}

\put(140,50){\makebox(0,0)[cc]{i=3j}}

\put(130,25){\makebox(0,0)[cc]{}}

\put(100,105){\makebox(0,0)[cc]{i=j}}

\put(125,55){\makebox(0,0)[cc]{}}

\linethickness{0.3mm}
\multiput(65,65)(0.12,-0.12){208}{\line(1,0){0.12}}
\put(90,40){\vector(1,-1){0.12}}
\linethickness{0.3mm}
\put(90,40){\line(0,1){20}}
\put(90,60){\vector(0,1){0.12}}
\linethickness{0.3mm}
\multiput(120,70)(0.12,-0.12){125}{\line(1,0){0.12}}
\put(135,55){\vector(1,-1){0.12}}
\linethickness{0.3mm}
\put(135,55){\line(0,1){20}}
\put(135,75){\vector(0,1){0.12}}
\linethickness{0.3mm}
\multiput(105,65)(0.12,-0.12){125}{\line(1,0){0.12}}
\put(120,50){\vector(1,-1){0.12}}
\linethickness{0.3mm}
\put(120,50){\line(0,1){20}}
\put(120,70){\vector(0,1){0.12}}
\linethickness{0.3mm}
\multiput(90,60)(0.12,-0.12){125}{\line(1,0){0.12}}
\put(105,45){\vector(1,-1){0.12}}
\linethickness{0.3mm}
\put(105,45){\line(0,1){20}}
\put(105,65){\vector(0,1){0.12}}
\end{picture}
\qed

To prove Proposition~\ref{estimates-worse},
we quote the following lemma in \cite{laff-delasalle}, which will be applied in the proof to some finite subgroups in $H.$

\begin{lem}\label{lem-single-pair}
(Lemma~4.9 in \cite{laff-delasalle}) 
Let $m,n\in\N^*,k\in\{1,2,...,m\}.$ Let $H$ be a locally compact amenable group, $\alpha,\beta:(\resr{m})^{n+1}\to H$ two maps. Let $f\in C_c(H)$ satisfying 
$$f(\alpha(a_1,...,a_n,b)\beta(x_1,...,x_n,y))=\lambda,
\tx{ if }y=\sum_{i=1}^n a_ix_i+b+\pi^k\in\resr{m},$$
and
$$f(\alpha(a_1,...,a_n,b)\beta(x_1,...,x_n,y))=\mu,\tx{ if }y=\sum_{i=1}^n a_ix_i+b+\pi^{k-1}\in\resr{m}.$$
Then 
$$|\lambda-\mu|\leq 2q^{-nk/2} \|f\|_{MA(H)}$$
where 
$$\|f\|_{MA(H)}=
\sup\{\|\sch{f}(\varphi)\|_{C_r^*(H)},\|\varphi\|_{C^*_r(H)}\leq 1\}.$$
\end{lem} 
Let us remark that
when $H$ is an arbitrary locally compact group, the lemma above still holds if in the conclusion $\|f\|_{MA(H)}$ is replaced by
$$\|f\|_{M_0A(H)}=
\sup_{B~H-C^*-\tx{algebra}}\{\|\sch{f}(\varphi)\|_{C_r^*(H,B)},\|\varphi\|_{C^*_r(H,B)}\leq 1\}.$$
But we do not need this generality in the proof of Proposition~\ref{estimates-worse}.

{\bf Proof of the first inequality in Proposition~\ref{estimates-worse}:}
We adapt the arguments in the proof of Lemma~2.1 in \cite{laff-preprint} to our situation by discretizing the matrices.

Denote $[\cdot]:\F_q((\pi))\to\F_q[\pi^{-1}]$ the integral part of an element as defined in the previous section.

Let $\sigma:\resr{i+1}\to\OO=\F_q[[\pi]]$ be any section. 
Define $\alpha,\beta:(\resr{i+1})^2\to H_1$ by 
$$\alpha(a,b)
=\begin{pmatrix}1&0&[\pi^{-i}\sigma(a)]&[\pi^{-i}\sigma(a^2-b)]
\\0&1&\pi^{-i}&[\pi^{-i}\sigma(a)]
\\0&0&1&0
\\0&0&0&1\end{pmatrix},$$
$$\beta(x,y)
=\begin{pmatrix}1&0&[\pi^{-i}\sigma(x/2)]&[\pi^{-i}\sigma(x^2/4+y)]
\\0&1&0&[\pi^{-i}\sigma(x/2)]
\\0&0&1&0
\\0&0&0&1\end{pmatrix},$$
for any $a,b,x,y\in\resr{i+1}.$

Compute $\alpha(a,b)\beta(x,y)$
$$=\begin{pmatrix}1&0&[\pi^{-i}(\sigma(x/2)+\sigma(a))]&[\pi^{-i}(\sigma(x^2/4+y)+\sigma(a^2-b))]
\\0&1&\pi^{-i}&[\pi^{-i}(\sigma(x/2)+\sigma(a))]
\\0&0&1&0
\\0&0&0&1\end{pmatrix}.$$

We see that $\forall a,b,x,y\in\resr{i+1},\|\alpha(a,b)\beta(x,y)\|=q^i$ (the $(2,3)$ matrix element achieves the maximal norm). Moreover, we see that
$$\det\begin{pmatrix}[\pi^{-i}(\sigma(x/2)+\sigma(a))]&[\pi^{-i}(\sigma(x^2/4+y)+\sigma(a^2-b))]
\\\pi^{-i}&[\pi^{-i}(\sigma(x/2)+\sigma(a))]\end{pmatrix}$$
$$=-\pi^{-2i}(y-ax-b)\tx{ mod }\pi^{-i+1}\OO.$$
So for any $l\in\{0,1,...,i\},$ when 
\begin{equation}\label{y-ax-b}
y-ax-b\in\pi^l\OO^\times/\pi^{i+1}\OO
\end{equation}
(where $\OO^\times$ denotes the group of unites of $\OO$), 
we have 
$$\|\wedge^2 (\alpha(a,b)\beta(x,y))\|=q^{2i-l}.$$
That is to say for any $l\in\{0,1,...,i\},$ when $a,b,x,y\in\resr{i+1}$ satisfy \ref{y-ax-b}, we have 
$$\alpha(a,b)\beta(x,y)\in KD(i,i-l)K.$$

Now denote for any $n\in\N^*,$
$$H_1^n=\{\begin{pmatrix}1 &0&x&z\\
 & 1 &y&x\\
 & & 1 &0\\
 & & & 1
\end{pmatrix},x,y,z\in\F_q[\pi^{-1}],|x|,|y|,|z|\leq q^n\},$$
which is a finite subgroup of $H_1.$ Note that the images of $\alpha$ and $\beta$ both lie in $H_1^i.$ Apply Lemma~\ref{lem-single-pair} to $n=1,m=i+1,k=i-j,$ 
$H=H_1^i,$ $\alpha,\beta$ as above, and $f=c|_{H_1^i},\lambda=c(D(i,j)),\mu=c(D(i,j+1)),$ we have
$$|c(D(i,j))-c(D(i,j+1))|\leq 2q^{-\frac{i-j}{2}}\|\sch{c|_{H_1^i}}\|
\leq\|\sch{c}|_{\C(H_1\cap B_{2i})}\|.$$
The last inequality is due to the facts that $H_1^i\subset H_1\cap B_{2i}$ and $\|f\|_{C^*_r(H_1^i)}=\|f\|_{C^*_r(H_1)},\forall f\in\C(H_1^i)\subset\C(H_1)$ (since $H_1$ is a discrete group).
\qed

{\bf Proof of the second inequality in Proposition~\ref{estimates-worse}:}
We will use discretization as in the proof of the first inequality, and improve the matrices in the proof of
Lemma~2.2 in \cite{laff-preprint}.
This improvement allows us to use the case of $n=2$ of Lemma~\ref{lem-single-pair} (whereas in the proof of the first inequality only $n=1$ can be used), 
resulting in the better factor $q^{-j}.$

We first write the matrices that are useful in both cases
when $i+j$ is even and when $i+j$ is odd.
Let $m\in\N.$
 Let $\sigma:\resr{m+1}\to\OO$ be any section. Define $\alpha,\beta:(\resr{m+1})^3\to H_2$ by 
$$\alpha(a,b)$$$$
=\begin{pmatrix}1&-[\pi^{-m-1}(1+\pi\sigma(a_1))]&[\pi^{-m-1}(1+\pi\sigma(a_2))]&-[\pi^{-2m}\sigma(b)]
\\0&1&0&[\pi^{-m-1}(1+\pi\sigma(a_2))]
\\0&0&1&[\pi^{-m-1}(1+\pi\sigma(a_1))]
\\0&0&0&1\end{pmatrix},$$
$$\beta(x,y)$$$$
=\begin{pmatrix}1&[\pi^{-m}\sigma(x_2)]&[\pi^{-m}\sigma(x_1)]&
\pi^{-m-1}[\pi^{-m}(\sigma(x_1)+\sigma(x_2))]+[\pi^{-2m}\sigma(y)]
\\0&1&0&[\pi^{-m}\sigma(x_1)]
\\0&0&1&-[\pi^{-m}\sigma(x_2)]
\\0&0&0&1\end{pmatrix},$$
for $a_1,a_2,b,x_1,x_2,y\in\resr{m+1}.$

Compute $$\alpha(a,b)\beta(x,y)
=\begin{pmatrix}1&-\xi_1&\xi_2&\eta
\\0&1&0&\xi_2
\\0&0&1&\xi_1
\\0&0&0&1\end{pmatrix},$$
where 
$$\xi_1=[\pi^{-m-1}(1+\pi\sigma(a_1)-\pi\sigma(x_2))],$$
$$\xi_2=[\pi^{-m-1}(1+\pi\sigma(a_2)+\pi\sigma(x_1))],$$
$$\eta=[\pi^{-2m}(\sigma(y)-\sigma(b))]-[\pi^{-m}\sigma(a_1)][\pi^{-m}\sigma(x_1)]
-[\pi^{-m}\sigma(a_2)][\pi^{-m}\sigma(x_2)]$$
$$=\pi^{-2m}(y-b-a_1x_1-a_2x_2)\tx{ mod }\pi^{-m+1}\OO.$$

Let us now prove the estimate when $i+j\in 2\N.$ Let $$m=(i+j)/2-1.$$
We see that for any $a_1,a_2,b,x_1,x_2,y\in\resr{m+1},$ 
$$\|\wedge^2\big(\alpha(a_1,a_2,b)\beta(x_1,x_2,y)\big)\|=q^{2m+2}=q^{i+j}.$$ 
Moreover, when 
\begin{equation}\label{y-a1x1-a2x2-b}
y-(a_1x_1+a_2x_2+b)\in\pi^{l}\OO^\times/\pi^{m+1}\OO,l\in\{0,1,...,m-1\},
\end{equation}
we have
$$\|\alpha(a_1,a_2,b)\beta(x_1,x_2,y)\|=q^{2m-l}.$$ 
In sum, for $a_1,a_2,b,x_1,x_2,y\in\resr{m+1}$ satisfying \eqref{y-a1x1-a2x2-b} above, we have 
$$\alpha(a_1,a_2,b)\beta(x_1,x_2,y)\in KD(i+j-2-l,l+2)K.$$ 

Denote for $n\in\N,$
$$H_2^n=\{\begin{pmatrix}1 &x&y&z\\
 & 1 &0&y\\
 & & 1 &-x\\
 & & & 1
\end{pmatrix},x,y,z\in\F_q[\pi^{-1}],|x|,|y|\leq q^n,|z|\leq q^{2n}\}.$$ It is a finite subgroup of $H_2.$ Note that the images of $\alpha$ and $\beta$ are both in $H_2^{(i+j)/2}.$ Apply Lemma~\ref{lem-single-pair} to $n=2,m=(i+j)/2,k=j-2,H=H_2^{(i+j)/2},$ and $\alpha,\beta$ as above, and $f=c|_{H_2^{(i+j)/2}},\lambda=c(D(i,j)),\mu=c(D(i+1,j-1)),$ and since $H_2^{(i+j)/2}\subset H_2\cap B_{i+j},$
we have 
$$|c(D(i,j))-c(D(i+1,j-1))|$$$$\leq 2q^{-(j-2)}
\|\sch{c|_{H_2^{(i+j)/2}}}\|$$ 
$$\leq 2q^2\cdot q^{-j}\|\sch{c}|_{\C(H_2\cap B_{i+j})}\|.$$

Now prove the estimate when $i+j\in 2\N+1.$ In this case let $$m=(i+j-1)/2-1.$$ 

Define the imbedding $\iota:\F_q\to H$ by 
$$\iota(\eps)=\begin{pmatrix}
1&0&0&0\\
&1&\pi^{-1}\eps&0\\
&&1&0\\
&&&1
\end{pmatrix},\forall\eps\in\F_q.$$
Define $H_2'^{n},n\in\N$ as the following subgroup of $H_2,$
$$H_2'^{n}=\{
\begin{pmatrix}
1&x&y&z\\
&1&0&y\\
&&1&-x\\
&&&1
\end{pmatrix}
|x|\leq q^n,|y|\leq, q^{n+1},|z|\leq q^{2n+1}\}.$$ 
$H_2$ is stable under the conjugate action of $\iota(\F_q),$ and so is $H_2'^n.$ Form $\tilde H_2=\iota(\F_q)\cdot
H_2,$ and $\tilde H_2^n=\iota(\F_q)\cdot
H_2'^n.$ $\tilde H_2$ is a subgroup of $H$ which contains $H_2$ as a subgroup of index $q.$

Now let $\tilde \alpha:(\resr{m+1})^3\to \tilde H_2$ be the map defined by 
$$\tilde \alpha(a_1,a_2,b)=\iota(1)\alpha(a_1,a_2,b),\forall a_1,a_2,b\in\resr{m+1}.$$ 
By easy computation we see that $\forall a_1,a_2,b,x_1,x_2,y\in\resr{m+1},$
$$\tilde\alpha(a_1,a_2,b)\beta(x_1,x_2,y)$$
$$=\begin{pmatrix}1&-\xi_1&\xi_2&\eta
\\0&1&\pi^{-1}&\xi_2+\pi^{-1}\xi_1
\\0&0&1&\xi_1
\\0&0&0&1\end{pmatrix},$$
We see 
$$\|\wedge^2 (\tilde \alpha(a_1,a_2,b)\beta(x_1,x_2,y))\|=q^{2m+3}=q^{i+j}.$$
And when 
$$y-(a_1x_1+a_2x_2+b)\in\pi^{l}\OO^\times/\pi^{m+1}\OO,l\in\{0,1,...,m-1\},$$
we have
$$\|\tilde \alpha(a_1,a_2,b)\beta(x_1,x_2,y)\|=q^{2m-l}.$$ Namely in this case, we obtain the decomposition
$$\tilde \alpha(a_1,a_2,b)\beta(x_1,x_2,y)\in KD(i+j-l-3,l+3)K.$$ 

The images of $\tilde \alpha$ and $\beta$ are both in $\tilde H_2^{(i+j-1)/2}.$ Now apply Lemma~\ref{lem-single-pair} to $n=2,m=(i+j-1)/2,k=j-3,H=\tilde H_2^{(i+j-1)/2},\tilde\alpha,\beta$ and $f=c|_{\tilde H_2^{(i+j-1)/2}},\lambda=c(D(i,j)),\mu=c(D(i+1,j-1)),$ and since $\tilde H_2^{(i+j-1)/2}\subset H_2\cap B_{i+j},$ we have 
$$|c(D(i,j))-c(D(i+1,j-1))|$$$$\leq 2q^{-(j-3)}
\|\sch{c|_{\tilde H_2^{(i+j)/2}}}\|$$ 
$$\leq 2q^{3}\cdot q^{-j}\|\sch{c}|_{\C(\tilde H_2\cap B_{i+j})}\|.$$
\qed

\section{Second proof of Theorem~\ref{main}}
In this section, another proof of Theorem~\ref{main}, more precisely of Theorem~\ref{obstruction}, is given by showing the following proposition (a slightly improved version of Proposition~\ref{estimates-worse} in the first inequality).

\begin{prop}\label{estimates}
For any $K$ biinvariant function $c\in C(G),$ we have
\begin{equation}\label{estimate-abelian}
|c(D(i,j))-c(D(i,j+1))|\leq 2q^{-(i-j)/2}\|\sch{c}|_{\C(H\cap B_{n_1})}\|,
\end{equation}
for any $(i,j)\in\Lambda$ with $i\geq 1$ 
where $$n_1=\max(\ell(D(i,j)),\ell(D(i,j+1)))=i+j+1$$ and
\begin{equation}\label{estimate-heisenberg}
|c(D(i,j))-c(D(i+1,j-1))|\leq 2q^2\cdot q^{-j}\|\sch{c}|_{\C(H\cap B_{n_2})}\|,
\end{equation}
for any $(i,j)\in\Lambda$ with $j\geq 3$ 
where $$n_2=\max(\ell(D(i,j)),\ell(D(i+1,j-1)))=i+j.$$
\end{prop}

Let $H_1$ be the abelian subgroup in $H$
$$H_1=\{h_1(x,y,z)=\begin{pmatrix}1 &0&x&z\\
 & 1 &y&x\\
 & & 1 &0\\
 & & & 1
\end{pmatrix},x,y,z\in\F_q[\pi^{-1}]\},$$ and let $H_2$ be the subgroup of Heisenberg type in $H$
$$H_2=\{h_2(x,y,z)=\begin{pmatrix}1 &x&y/2&z\\
 & 1 &0&y/2\\
 & & 1 &-x\\
 & & & 1
\end{pmatrix},x,y,z\in\F_q[\pi^{-1}]\}.$$ 
The group law is as follows:
$$h_2(x,y,z)h_2(x',y',z')=h_2(x+x',y+y',z+z'+\frac{1}{2}(xy'-yx')).$$

The proof of the two inequalities in Proposition \ref{estimates} replies on the construction as follows of two family of explicit functions on $H_1$ and $H_2$ respectively.

Denote $[\cdot]:\F_q((\pi))\to\F_q[\pi^{-1}]$ the integral part of an element, i.e.
$$[a_i\pi^{-i}+a_{i-1}\pi^{-i+1}+...+a_1\pi^{-1}+a_0+a_{-1}\pi+...]$$
$$=a_i\pi^{-i}+a_{i-1}\pi^{-i+1}+...+a_1\pi+a_0,\forall a_*\in\F_q.$$
Now fix $(i,j)\in\Lambda.$
Define
$$h_{1,i,j}=\esp{a\in\resr{i}}e_{h_1([\pi^{-i}a],\pi^{-i},[\pi^{-i}a^2]+\pi^{-j})},$$
$$h_{2,i,j}=\esp{a,b,c\in\resr{i}}e_{h_2([\pi^{-m}(1+\pi a)],[\pi^{-m}b],[\pi^{-i}(1+\pi c)])},$$
where $m=m_{i,j}$ is the integral part of $(i+j)/2,$ i.e. when $i+j\in 2\N,m=(i+j)/2,$ when $i+j\in 2\N+1, m= (i+j-1)/2.$
More explicitly,
$$h_1([\pi^{-i}a],\pi^{-i},[\pi^{-i}a^2]+\pi^{-j})$$$$
=\begin{pmatrix}1 &0&[\pi^{-i}a]&[\pi^{-i}a^2]+\pi^{-j}\\
 & 1 &\pi^{-i}&[\pi^{-i}a]\\
 & & 1 &0\\
 & & & 1
\end{pmatrix},a\in\resr{i}.$$
$$h_2([\pi^{-m}(1+\pi a)],[\pi^{-m}b],[\pi^{-i}(1+\pi c)])$$$$
=\begin{pmatrix}1 &[\pi^{-m}(1+\pi a)]&[\pi^{-m}b]/2&[\pi^{-i}(1+\pi c)]\\
 & 1 &0&[\pi^{-m}b]/2\\
 & & 1 &-[\pi^{-m}(1+\pi a)]\\
 & & & 1
\end{pmatrix},a,b,c\in\resr{i}.$$
The explicit functions are defined as
\begin{align}
\Delta_{1,i,j}=h_{1,i,j}-h_{1,i,j+1}\in\C H_1,\label{explicit-functions-1}\\
\Delta_{2,i,j}=h_{2,i,j}-h_{2,i+1,j-1}\in\C H_2. \label{explicit-functions-2}
\end{align} 


\begin{prop}\label{explicit}
\begin{equation}\label{explicit-abelian}
\|\Delta_{1,i,j}\|_{C_r^*(H_1)}\leq 2q^{-(i-j)/2}
\end{equation}
\begin{equation}\label{explicit-heisenberg}
\|\Delta_{2,i,j}\|_{C_r^*(H_2)}\leq 2q^2\cdot q^{-j}
\end{equation}
\end{prop}

{\bf Proof of Proposition~\ref{estimates} by Proposition~\ref{explicit}:}
Recall that for any $g=(g_{\alpha,\beta})_{1\leq \alpha,\beta\leq 4}\in G,$ $g\in KD(i,j)K$ for $(i,j)\in\Lambda$ if and only if 
$$\|g\|=\max_{1\leq \alpha,\beta\leq 4}|g_{\alpha,\beta}|=q^{i}$$ and 
$$\|\wedge^2 g\|=\max_{1\leq \alpha_1,\beta_1,\alpha_2,\beta_2\leq 4}|g_{\alpha_1,\beta_1}g_{\alpha_2,\beta_2}-g_{\alpha_1,\beta_2}g_{\alpha_2,\beta_1}|=q^{i+j}.$$

From the construction \eqref{explicit-functions-1} we see that
$$h_1([\pi^{-i}a],\pi^{-i},[\pi^{-i}(a^2+\pi^{i-j})])\in H_1\cap KD(i,j)K,$$
i.e.
$$\supp h_{1,i,j}\subset H_1\cap KD(i,j)K.$$

Since $H$ is amenable, we have
$$|c(D(i,j))-c(D(i,j+1))|=|\sum_{h\in H}c(h)\Delta_{1,i,j}(h)|$$
$$\leq \|\sch{c}(\Delta_{1,i,j})\|_{C_r^*(H)}
\leq\|\sch{c}|_{\C(H\cap B_{n_1})}\|\|\Delta_{1,i,j}\|_{C^*_r(H)}.$$

Now that $H_1$ is a subgroup of $H,$ we have that $\|\Delta_{1,i,j}\|_{C_r^*(H)}=\|\Delta_{1,i,j}\|_{C_r^*(H_1)},$ so the first inequality is proved.

The second inequality requires a bit more computations.

First when $i+j\in 2\N,$ we have (see \eqref{explicit-functions-2})
$$\supp h_{2,i,j}\subset H_2\cap KD(i,j)K,$$

$$|c(D(i,j))-c(D(i+1,j-1))|=|\sum_{h\in H}c(h)\Delta_{2,i,j}(h)|$$
$$\leq \|\sch{c}(\Delta_{2,i,j})\|_{C_r^*(H)}
\leq\|\sch{c}|_{\C(H\cap B_{n_2})}\|\|\Delta_{2,i,j}\|_{C^*_r(H)},$$
and
$$\|\Delta_{2,i,j}\|_{C_r^*(H)}=\|\Delta_{2,i,j}\|_{C_r^*(H_2)}.$$

When $i+j\in 2\N+1,$
$$\iota(1)h_2([\pi^{-m}(1+\pi a)],[\pi^{-m}b],[\pi^{-i}(1+\pi c)])$$$$
=\begin{pmatrix}1 &[\pi^{-m}(1+\pi a)]&[\pi^{-m}b]/2&[\pi^{-i}(1+\pi c)]\\
 & 1 &\pi^{-1}&[\pi^{-m}b]/2-\pi^{-1}[\pi^{-m}(1+\pi a)]\\
 & & 1 &-[\pi^{-m}(1+\pi a)]\\
 & & & 1
\end{pmatrix}\in KD(i,j)K,$$
$\forall a,b,c\in\resr{i},$ where $m=(i+j-1)/2$ as before, and
$$\iota(\eps)
=\begin{pmatrix}1 &0&0&0\\
 & 1 &\eps\pi^{-1}&0\\
 & & 1 &0\\
 & & & 1
\end{pmatrix},\eps\in\F_q.$$

Finally we have
$$|c(D(i,j))-c(D(i+1,j-1))|=|\sum_{h\in H}c(h)\Delta_{2,i,j}(\iota(1)h)|$$
$$\leq \|\sch{c}(L_{\iota(-1)}\Delta_{2,i,j})\|_{C_r^*(H)}
\leq\|\sch{c}|_{\C(H\cap B_{n_2})}\|\|L_{\iota(-1)}\Delta_{2,i,j}\|_{C^*_r(H)},$$
and
$$\|L_{\iota(-1)}\Delta_{2,i,j}\|_{C_r^*(H)}=\|\Delta_{2,i,j}\|_{C_r^*(H_2)},$$
so the second inequality is proved.

\qed

Now it suffices to show Proposition \ref{explicit}, whose proof unlike Proposition \ref{estimates-worse} does not rely on Lemma 4.9 in \cite{laff-delasalle}.

{\bf Proof of inequality \eqref{explicit-abelian} in Proposition~\ref{explicit}:}
Here we follow \cite{laff-orsay}.

\begin{lem}(norm of quadratic Gauss sum)\label{gauss-sum}
If character $\eta\in\widehat{\resr{\ell}}$ is non-degenerate 
(i.e. $\eta|_{\pi^{\ell-1}\OO/\pi^\ell\OO}\neq 1$),
then 
$$|\esp{a\in\resr{\ell}}\eta(a^2)|=q^{-\ell/2}.$$
\end{lem}
{\bf Proof}
$$|\esp{a\in\resr{\ell}}\eta(a^2)|^2
=|\esp{a,b\in\resr{\ell}}\eta(a^2-b^2)|
=|\esp{a,b\in\resr{\ell}}\eta\big((a-b)(a+b)\big)|.$$
Since $q$ is odd, we can introduce new variables $x=a+b,y=a-b$ which is an invertable linear transform on $(\resr{\ell})^2,$ thus
$$|\esp{a\in\resr{\ell}}\eta(a^2)|^2=|\esp{x,y\in\resr{\ell}}\eta(xy)|=q^{-\ell}.$$ 
\qed

Since $H_1$ is an abelian group,
$\forall\varphi\in\C H_1,$
$$\|\varphi\|_{C_r^*(H_1)}=\sup_{\chi\in \widehat{H_1}}
\big|\chi(\varphi)\big|.$$

Fix $\chi\in \widehat{H_1},$ and suppose $\chi_1,\chi_2,\chi_3\in\widehat{\F_q[\pi^{-1}]}$ such that
$\forall x,y,z\in\F_q[\pi^{-1}],\chi(h_1(x,y,z))=\chi_1(x)\chi_2(y)\chi_3(z).$

$$\chi(\Delta_{1,i,j})=\chi_2(\pi^{-i})\Big(\chi_3(\pi^{-j})-\chi_3(\pi^{-j-1})\Big)
\esp{a\in\resr{i}}\chi_1([\pi^{-i}a])\chi_3([\pi^{-i}a^2]).$$

We see that unless $\chi(\Delta_{1,i,j})=0,$ we have 
$\Ker\big(\chi_3([\pi^{-i}\cdot])\big)\subset \pi^{i-j}\OO,$
and
$\Ker(\chi_1([\pi^{-i}\cdot]))\supset\Ker(\chi_3([\pi^{-i}\cdot])).$
\footnote{
If we replace $h_{1,i,j}$ by the function
$$h'_{1,i,j}=\esp{a,b,c\in\resr{i}}
e_{h_1([\pi^{-i}a],[\pi^{-i}(1+\pi b)],[\pi^{-i}a^2+\pi^{-j}(1+\pi c)])},$$
we can then locate the support of $\Delta'_{1,i,j}=h'_{1,i,j}-h'_{1,i,j+1}$ more precisely, i.e. we have
$$\Ker\big(\chi_3([\pi^{-i}\cdot])\big)=\pi^{i-j}\OO \tx{ or } \pi^{i-j+1}\OO,$$
$$\Ker\big(\chi_1([\pi^{-i}\cdot])\big)\supset\pi^{i-j+1}\OO$$
and
$$\Ker(\chi_2([\pi^{-i}\cdot]))\supset\pi\OO.$$
But this is not very useful in the current situation.}
Consequently, there exists $\theta\in\OO,$ such that
$\chi_1([\pi^{-i}\cdot])=\chi_3([\pi^{-i}\theta\cdot]).$

Now we have
$$|\chi(\Delta_{1,i,j})|\leq 2\Big|\esp{a\in\resr{i}}\chi_3([\pi^{-i}(\theta a+a^2)])\Big|.$$

Since $q$ is odd,
$(\theta/2)^2+\theta a+a^2=(\theta/2+a)^2,|\chi_3([\pi^{-i}(\theta/2)^2])|=1,$
we have by Lemma \ref{gauss-sum},
$$|\chi(\Delta_{1,i,j})|\leq 2\Big|\esp{a\in\resr{i}}\chi_3([\pi^{-i}(\theta/2+a)^2])\Big|\leq 2q^{-(i-j)/2}.
\footnote{
In fact, it suffices to prove this final inequality for rational $\chi\in\widehat{H_1},$ 
i.e. there exist $\psi\in\widehat{F}$ vanishing on $\pi\OO$ and being non zero on $\OO,$ 
and $t_1,t_2,t_3\in\F_q[\pi^{-1}]$ such that 
$\chi_i(\cdot)=\psi(t_i\cdot),i=1,2,3.$
Note that $\psi(t\cdot)$ vanishes on $\pi^{\ell+1}\OO$ and is non zero on $\pi^{\ell}\OO$
if and only if $|t|=q^{\ell}.$
With this notation, $|t_3|\geq q^{j-1},|t_1|\leq |t_3|$ unless $\chi(\Delta_{1,i,j})=0.$
}
$$
\qed

{\bf Proof of inequality \eqref{explicit-heisenberg} in Proposition~\ref{explicit}:}

Let $\chi,\chi'\in\widehat{\F_q[\pi^{-1}]},\chi\neq 0.$ 
We define a unitary representation of $\rho_{\chi,\chi'}:H_2\to\mathcal U(\ell^2(\F_q[\pi^{-1}]))$
by
$$\rho_{\chi,\chi'}(h_2(a,b,c))f(x)=f(x+a)\chi(xb)\chi(c+\frac{1}{2}ab)\chi'(b)$$
(it is well defined since $q$ is odd).

V.~Lafforgue suggested the following formula for calculating the $C_r^*$ norms on $H_2.$

\begin{lem}\label{lem-plancherel}
$\forall \varphi\in\C H_2,$
$$\|\varphi\|_{C^*_r(H_2)}=\sup_{\chi,\chi'\in\widehat{\F_q[\pi^{-1}]},\chi\neq 0}
\|\rho_{\chi,\chi'}(\varphi)\|_{\mathcal L(\ell^2(\F_q[\pi^{-1}]))}.$$
\end{lem}
{\bf Remarks.}\begin{enumerate}
\item Being a counterpart of $H_2$ in a number field,
the following discrete Heisenberg group
$$H_2(\Z)=\{h_2(x,y,z)=\begin{pmatrix}1 &x&y/2&z\\
 & 1 &0&y/2\\
 & & 1 &-x\\
 & & & 1
\end{pmatrix},x,z\in \Z, y\in 2\Z\}$$
also admits a similar formula for the $C^*_r$ norms. More precisely,
for $\theta,\theta' \in[0,1),$
define the unitary representation $\rho_{\theta,\theta'}:H_2(\Z)\to\mathcal U(\ell^2(\Z))$
of central character $\theta$
by a similar formula
 $$\rho_{\theta,\theta'}(h_2(a,b,c))f(x)
=f(x+a)e^{2\pi i \theta xb}e^{2\pi i\theta (c+\frac{1}{2}ab)}e^{2\pi i\theta'b},$$
then we have $\forall \varphi\in\C (H_2(\Z)),$
$$\|\varphi\|_{C^*_r(H_2(\Z))}=\sup_{\theta,\theta'\in[0,1)}
\|\rho_{\theta,\theta'}(\varphi)\|_{\mathcal L(\ell^2(\Z))}.$$
\item The formula for the $C^*_r$ norm in the first remark can be reduced 
to those irrational $\theta\in[0,1)\backslash\Q$ and $\theta'=0,$ i.e.
$\forall \varphi\in\C (H_2(\Z)),$
$$\|\varphi\|_{C^*_r(H_2(\Z))}=\sup_{\theta\in[0,1)\backslash\Q}
\|\rho_{\theta,0}(\varphi)\|_{\mathcal L(\ell^2(\Z))}$$
(The analogues formula also holds for $H_2$, but we don't use it the proof).
Indeed, 
when $\theta$ is irrational, $\rho_{\theta,0}(H_2(\Z))$ generates algebra $A_\theta$ 
which is a simple $C^*$ algebra, i.e. any representation of $A_\theta$ is faithful.
Moreover, for any $C^*$ algebra $A$ 
and any representation $\sigma_1,\sigma_2:A\to\mathcal L(H),$ we have
$$\Ker\sigma_1\subset\Ker\sigma_2\Leftrightarrow 
\|\sigma_1(a)\|_{\mathcal L(H)} \geq\|\sigma_2(a)\|_{\mathcal L(H)},
\forall a\in A.$$
By applying the previous fact to the representation of $A_\theta$ generated by  $\rho_{\theta,\theta'}(H_2(\Z)),$ we have
$$\|\rho_{\theta,\theta'}(\varphi)\|=\|\rho_{\theta,0}(\varphi)\|,\forall\theta'\in[0,1),
\forall \varphi\in\C(H_2(\Z)).$$
\end{enumerate}

{\bf Proof of the lemma.} Let $N_2\supset H_2$ be the following Heisenberg group
$$N_2=\{h_2(x,y,z)=\begin{pmatrix}1 &x&y/2&z\\
 & 1 &0&y/2\\
 & & 1 &-x\\
 & & & 1
\end{pmatrix},x,y,z\in F=\F_q((\pi))\},$$ 
and for a character $\eta\in\hat F\backslash\{0\},$ denote $\rho_\eta:N_2\to\mathcal U(L^2(F))$ the representation defined by
$$\rho_{\eta}(h_2(a,b,c))f(x)=f(x+a)\eta(xb)\eta(c+\frac{1}{2}ab),a,b,c,x\in F.$$

Let $D$ be the fundamental domain $D=\pi\OO$ for the translation of $\F_q[\pi^{-1}]$ on $F.$

We have an isomorphism of representations of $H$
$$\rho_\eta|_{H_2}\simeq\int_D^\oplus \rho_{\eta|_{\F_q[\pi^{-1}]},\eta|_{\F_q[\pi^{-1}]}(\delta\cdot)}d\delta$$
defined by
$$L^2(F)\xrightarrow{\sim}\int_D^\oplus\ell^2(\F_q[\pi^{-1}])d\delta,$$
$$\phi\mapsto\Big(\big[r\mapsto\phi(r+\delta)\big]\in\ell^2(\F_q[\pi^{-1}])\Big)_{\delta\in D},$$
where $\eta|_{\F_q[\pi^{-1}]}(\delta\cdot)$ denotes the character $[\gamma\mapsto\eta(\delta\gamma)]\in\widehat{\F_q[\pi^{-1}]}.$
\qed

We write the action of $\Delta_{2,i,j}$ in the following form
$$\rho_{\chi,\chi'}(\Delta_{2,i,j})f(x)$$
$$=C \esp{a,b\in\OO}\Big(f(x+[\pi^{-m}(1+\pi a)])
\chi\big((x+\frac{1}{2}[\pi^{-m}(1+\pi a)])[\pi^{-m}b]\big)\chi'([\pi^{-m}b])\Big),$$
where
$$C=\esp{c\in\OO}\Big(\chi([\pi^{-i}(1+\pi c)])-\chi([\pi^{-i-1}(1+\pi c)])\Big).$$

We use excessively the following basic fact in the proof: 
for any finite abelian group $A$ and any unitary character $\eta\in\widehat{A},$ we have
\begin{eqnarray}\label{char}
\esp{a\in A}\eta(a)=1 \tx{ when } \eta\in\widehat{A} \tx{ is trivial;}\nonumber\\
\esp{a\in A}\eta(a)=0 \tx{ when } \eta\in\widehat{A} \tx{ is non-trivial.}
\end{eqnarray}

\begin{lem}\label{kernel}
If $C\neq 0,$ then 
$\chi|_{[\pi^{-i+1}\OO]}=1\in\widehat{[\pi^{-i+1}\OO]}$ and 
$\chi|_{[\pi^{-i-1}\OO]}\neq 1\in\widehat{[\pi^{-i-1}\OO]}.$
\end{lem}
{\bf Proof of in Lemma~\ref{kernel}.}
If $\chi|_{[\pi^{-i-1}\OO]}$ is trivial, then 
$C=0$
 since
$\chi(z_1)-\chi(z_2)=0,\forall z_1\in[\pi^{-i-1}\OO],z_2\in[\pi^{-i}\OO].$
On the other hand, if $\chi|_{[\pi^{-i+1}\OO]}$ is non-trivial, then by \eqref{char}
$$\esp{c\in\OO}\chi([\pi^{-i}(1+\pi c)])=\chi(\pi^{-i})\esp{z\in[\pi^{-i+1}\OO]}\chi(z)=0,$$
$$\esp{c\in\OO}\chi([\pi^{-i-1}(1+\pi c)])\chi(\pi^{-i-1})\esp{z\in[\pi^{-i}\OO]}\chi(z)=0,$$
and therefore $C=0.$\qed

The matrix of $\rho_{\chi,\chi'}(\Delta_{2,i,j})$ is block diagonal 
and each block corresponds to a coset 
$x_0+[\pi^{-m}\OO],x_0\in\F_q[\pi^{-1}].$ 
Indeed, the action of $\rho_{\chi,\chi'}(\Delta_{2,i,j})$ on $\ell^2(\F_q[\pi^{-1}])$ only concerns translations of elements in $[\pi^{-m}\OO]$ and scalars. 

It remains to show that each block of $\rho_{\chi,\chi'}(\Delta_{2,i,j})$ associated to the coset $x_0+[\pi^{-m}\OO]$ has norm $\leq 2q^{2-j},$
$$\|\rho_{\chi,\chi'}(\Delta_{2,i,j})|_{\ell^2(x_0+[\pi^{-m}\OO])}\|_{\mathcal L(\ell^2(x_0+[\pi^{-m}\OO]))}
\leq 2q^{2-j}.~~~~~~~~~(*)$$ 

Now fix a coset $x_0+[\pi^{-m}\OO]$ for some $x_0\in\F_q[\pi^{-1}].$ 
We provide two proofs of $(*)$. The two proofs are related, but the author thinks that both have merits and it might be useful to write them down.

\noindent {\bf First proof of $(*)$}:

Denote $E_{\eps}$
the subset $x_0+\pi^{-m}\eps+[\pi^{-m+1}\OO]\subset x_0+[\pi^{-m}\OO]$
for $\eps\in\F_q.$
We have a disjoint union decomposition 
$$x_0+[\pi^{-m}\OO]=\sqcup_{\eps\in\F_q}E_{\eps}.$$
For each $\eps\in\F_q,$ 
$\rho_{\chi,\chi'}(\Delta_{2,i,j})$ sends $\ell^2(E_{\eps})$ to $\ell^2(E_{\eps-1}),$ 
and thus the action of $\rho_{\chi,\chi'}(\Delta_{2,i,j})$ on $\ell^2(x_0+[\pi^{-m}\OO])$ has the following form of block matrix
$$\begin{pmatrix}
0&*&0&...&0\\
0&0&*&...&0\\
&&...&&\\
0&0&0&...&*\\
*&0&0&...&0\\
\end{pmatrix},$$
where each block $*$ has size $q^{m-1}$ 
and corresponds to the action $\ell^2(E_{\eps})\to \ell^2(E_{\eps-1}).$

The following lemma claims that after appropriate identification of $E_{\eps}$ and $E_{\eps-1}$
the block $*$ corresponding to $\eps$ is 
$Cq^{-2m+1+i}(\simeq q^{-j})$ times 
the projection from $\ell^2(E_{\eps})$ onto $[\pi^{m-i}\OO]$ invariant vectors in $\ell^2(E_{\eps}),$
and thus our inequality follows.
More precisely, by identifying $x_0+\pi^{-m}(\eps-1)-y+y_{\eps}\in E_{\eps-1}$ and 
$x_0+\pi^{-m}\eps+y\in E_{\eps},$ 
$\rho_{\chi,\chi'}(\Delta_{2,i,j})$ sends $\delta_{x_0+\pi^{-m}\eps+y}$ to
$Cq^{-2m+1+i}\esp{\alpha\in\OO}\delta_{x_0+\pi^{-m}\eps+y+[\pi^{m-i}\alpha]},$
and thus has norm less than $ 2q^{-2m+1+i}\leq 2q^{2-j}.$

{\bf Remark.} The identification of $E_{\eps-1}$ and $E_{\eps}$ via
$$x_0+\pi^{-m}(\eps-1)+y_{\eps} -y \to
x_0+\pi^{-m}\eps+y $$
corresponds to the fact that $A_{x,y}$ is an anti-diagonal in the second proof below (the center of the anti-diagonal is $x_0+\pi^{-m}(\eps-\frac{1}{2})+\frac{1}{2}y_{\eps} $).

\begin{lem}
If $\rho_{\chi,\chi'}(\Delta_{2,i,j})|_{\ell^2(E_{\eps})}
\neq 0\in\mathcal L(\ell^2(E_{\eps}),\ell^2(E_{\eps-1})),$
then there exists $y_{\eps}\in [\pi^{-m+1}\OO],$ such that
$\forall y\in[\pi^{-m+1}\OO]$
$$\rho_{\chi,\chi'}(\Delta_{2,i,j})f(x_0+\pi^{-m}(\eps-1)+y)$$
$$=Cq^{-2m+1+i}\esp{\alpha\in\OO}f(x_0+\pi^{-m}\eps-y+y_{\eps}+[\pi^{m-i}\alpha]),$$
for any $ i\geq j\geq 2.$
\end{lem}

{\bf Proof of the lemma:}
By hypothesis there exist $f_0\in\ell^2(E_{\eps})$ and $y_0\in[\pi^{-m+1}\OO]$ such that
$$\rho_{\chi,\chi'}(\Delta_{2,i,j})f_0(x_0+\pi^{-m}(\eps-1)+y_0)$$
$$=C \esp{a,b\in\OO}\Big(f_0(x_0+\pi^{-m}\eps+y_0+[\pi^{-m+1} a)])$$
$$\chi\big((x_0+\pi^{-m}(\eps-1)+y_0+\frac{1}{2}[\pi^{-m}(1+\pi a)])[\pi^{-m}b]\big)\chi'([\pi^{-m}b])\Big)$$
$$\neq 0.$$
By fixing $a$ and averaging over $b$ we see that there exists $a_0\in\OO$ such that 
\begin{equation}\label{chi=1}
\chi\big((x_0+\pi^{-m}(\eps-1)+y_0+\frac{1}{2}[\pi^{-m}(1+\pi a_0)])[\pi^{-m}b]\big)\chi'([\pi^{-m}b])=1,
\forall b\in\OO.
\end{equation}
Set $y_{\eps}=[\pi^{-m+1}a_0]+2y_0\in[\pi^{-m+1}\OO].$

By definition $\forall f\in\ell^2(E_{\eps}),y\in[\pi^{-m+1}\OO]$ we have
$$\rho_{\chi,\chi'}(\Delta_{2,i,j})f(x_0+\pi^{-m}(\eps-1)+y)$$
$$=C \esp{a,b\in\OO}\Big(f(x_0+\pi^{-m}\eps+y+[\pi^{-m+1}a])$$
$$\chi\big((x_0+\pi^{-m}(\eps-1)+y+\frac{1}{2}[\pi^{-m}(1+\pi a)])[\pi^{-m}b]\big)\chi'([\pi^{-m}b])\Big),$$
by equality \eqref{chi=1} it equals
$$=C \esp{a,b\in\OO}\Big(f(x_0+\pi^{-m}\eps+y+[\pi^{-m+1}a])$$
$$\chi\big((y-y_0+\frac{1}{2}[\pi^{-m+1}(a-a_0)])[\pi^{-m}b]\big)\Big),$$
by change of variables $a'=a-a_0+2\pi^{m-1}(y-y_0)$ (where $2\pi^{m-1}(y-y_0)\in\F_q$) it equals
$$=C \esp{a',b\in\OO}\Big(f(x_0+\pi^{-m}\eps-y+y_{\eps}+[\pi^{-m+1}a'])
\chi\big(\frac{1}{2}[\pi^{-m+1}a'][\pi^{-m}b]\big)\Big).$$

When $\chi|_{[\pi^{-i}\OO]}=1,\chi|_{[\pi^{-i-1}\OO]}\neq 1,$
we have $\chi([\pi^{-m}c]/2)=1,\forall c\in\OO$ and then
$$\chi(\frac{1}{2}[\pi^{-m+1}a'][\pi^{-m}b])=\chi(\frac{1}{2}[\pi^{-2m+1}a'b]),\forall a',b\in\OO.$$
Thus
$$\rho_{\chi,\chi'}(\Delta_{2,i,j})f(x_0+\pi^{-m}(\eps-1)+y)$$
$$=C \esp{a',b\in\OO}\Big(f(x_0+\pi^{-m}\eps-y+y_{\eps}+[\pi^{-m+1}a'])
\chi\big(\frac{1}{2}[\pi^{-2m+1}a'b]\big)\Big).$$
Being a Fourier transform for the non-degenerate character 
$[\alpha\mapsto\chi(\frac{1}{2}[\pi^{-2m+1}\alpha])]\in\widehat{\resr{2m-1-i}},$ 
it equals
$$Cq^{-2m+1+i}\esp{\alpha\in\OO}f(x_0+\pi^{-m}\eps-y+y_{\eps}+[\pi^{m-i}\alpha]).$$

The case 
when $\chi|_{[\pi^{-i+1}\OO]}=1,\chi|_{[\pi^{-i}\OO]}\neq 1$
can be handled similarly.
\qed

This ends the first proof of $(*).$

\noindent {\bf Second proof of $(*)$} (due to V.~Lafforgue):
\begin{lem}\label{average}
If $\chi|_{[\pi^{-i+1}\OO]}=1\in\widehat{[\pi^{-i+1}\OO]}$ and 
$\chi|_{[\pi^{-i-1}\OO]}\neq 1\in\widehat{[\pi^{-i-1}\OO]},$ then
unless $w\in[\pi^{-(i-m)}\OO],$ we have that
$\forall w\in[\pi^{-m}\OO],$ 
$$\esp{z\in[\pi^{-m}\OO]}\chi(\frac{1}{2}wz)=0,$$
or equivalently by \eqref{char}
$$[z\mapsto\chi(\frac{1}{2}wz)]\in\widehat{[\pi^{-m}\OO]}\tx{ is non-trivial}.$$
\end{lem}
{\bf Proof of in Lemma~\ref{average}.}
We prove it by contradiction.
Suppose $w=\pi^{m-i-\alpha}w_0\in[\pi^{m-i-\alpha}\OO^\times],
w_0\in \F_q+...+\pi^{-m+i+\alpha}\F_q,
\alpha\in\{1,2,...,2m-i\},$ such that
$$\chi|_{w[\pi^{-m}\OO]}=1\in\widehat{[\pi^{-m}\OO]}.$$

We have 
$$\chi|_{\pi^{-i-1}w_0(\F_q+\F_q\pi)}=1.$$ Indeed,
since $1\leq \alpha\leq m,$ we have
$$\pi^{-i-1}w_0(\F_q+\F_q\pi)=\pi^{-i-\alpha}w_0(\F_q\pi^{\alpha-1}+\F_q\pi^\alpha)$$
$$\subset \pi^{-i-\alpha}w_0(\F_q+\F_q\pi+...+\F_q\pi^m)=w[\pi^{-m}\OO].$$

Now $\forall\eps_1,\eps_2\in\F_q,$ there exist $\eps_1',\eps_2'\in\F_q$ such that 
$\eps_1+\eps_2\pi\in w_0(\eps_1'+\eps_2'\pi)+\pi^2\OO.$ Since $\chi|_{[\pi^{-i+1}\OO]}=1,$
we have 
$$\chi(\pi^{-i-1}(\eps_1+\eps_2\pi))=\chi(\pi^{-i-1}w_0(\eps_1'+\eps_2'\pi))=1.$$
As a consequence $\chi|_{[\pi^{-i-1}\OO]}=1,$ which is a contradiction to 
the hypothesis in the lemma.
\qed

Let $A=(A_{x,y})_{x,y\in x_0+[\pi^{-m}\OO]}$ be the matrix
of the block of $\rho_{\chi,\chi'}(\Delta_{2,i,j})$ associated to $\ell^2(x_0+[\pi^{-m}\OO]).$ 

We will show that $\|A\|_{\mathcal L(\ell^2(x_0+[\pi^{-m}\OO]))}\leq 2q^{1-j}.$

First we have $A_{x,y}=0$ unless $y\in x+\pi^{-m}+[\pi^{-m+1}\OO],$ and in this case, (since $|[\pi^{-m}\OO]|=q^{m+1}$)
$$A_{x,y}=Cq^{-m-1}\esp{z\in[\pi^{-m}\OO]}\chi(\frac{x+y}{2}z)\chi'(z),$$ i.e. by \eqref{char}
$$A_{x,y}=Cq^{-m-1}\tx{ when }[z\mapsto\chi(\frac{x+y}{2}z)\chi'(z)]\in\widehat{[\pi^{-m}\OO]}\tx{ is trivial, and}$$
$$A_{x,y}=0\tx{ when }[z\mapsto\chi(\frac{x+y}{2}z)\chi'(z)]\in\widehat{[\pi^{-m}\OO]}\tx{ is non-trivial.}$$

Now suppose $x,y,y'\in x_0+[\pi^{-m}\OO]$ such that both $A_{x,y}$ and $A_{x,y'}$ are non-zero. 
By taking ratio we see that $[z\mapsto\chi(\frac{1}{2}(y-y')z)]\in[\pi^{-m}\OO]$ is a trivial character.
By Lemma~\ref{average} we see that $y-y'\in [\pi^{-(i-m)}\OO].$

By the same argument for $x,x',y\in x_0+[\pi^{-m}\OO]$ such that both $A_{x,y}$ and $A_{x',y}$ are non-zero,
we have $x-x'\in [\pi^{-(i-m)}\OO].$

Therefore, each line and column in $A$ has at most $|[\pi^{-(i-m)}\OO]|=q^{i-m+1}$ non-zero coefficients.
Each coefficient in $A$ has absolute value at most $2q^{-m-1}.$
The $\ell^2$ norm of $A$ is at most $2q^{-m-1}\cdot q^{i-m+1}=2q^{i-2m}\leq 2q^{-j+1},$ 
and so is the operator norm of $A.$

{\bf Remark 1.} By the same argument, for $x,x',y,y'\in x_0+[\pi^{-m}\OO]$ such that both $A_{x,y}$ and $A_{x',y'}$ are non-zero, we have $(x+y)-(x'+y')\in [\pi^{-(i-m)}\OO].$
It means that $A$ is a block ''anti-diagonal''.

{\bf Remark 2.} Following the previous remark, we can write the action of $A$ in the following form (supposing $A_{x,y}\neq 0$)
$$Af(x')=\sum_{y'\in x_0+[\pi^{-m}\OO]}A_{x',y'}f(y')
=\sum_{\alpha \in \OO}A_{x',x+y-x'+[\pi^{m-i}\alpha]}f(x+y-x'+[\pi^{m-i}\alpha]),$$
where $A_{x',x+y-x'+[\pi^{m-i}\alpha]}=0$ or $Cq^{-m-1},$ which means that $A$ is 
roughly (the precise formula requires a more detailed analysis on Lemma~\ref{average})
$Cq^{i-2m}$ 
times
the projection onto $[\pi^{m-i}\OO]$-invariant functions in $\ell^2(x_0+[\pi^{m-i}\OO])$,
after identifying $x'$ to $x+y-x',$
corresponding to the arguments in the first proof above.
\qed

\end{document}